\documentclass[11pt]{amsart}

\usepackage{amsmath}
\usepackage{amsthm}
\usepackage{amssymb}
\usepackage{mathrsfs}
\usepackage{url}

\renewcommand{\today}{\number\day\nobreakspace\ifcase\month\or
  January\or February\or March\or April\or May\or June\or
  July\or August\or September\or October\or November\or December\fi,
  \number\year}
\let\mathnoun\mathrm
\let\ptref\eqref

\newcommand{\hyp}{\text{-}}
\newcommand{\setm}{\setminus}
\newcommand{\then}{\rightarrow}
\newcommand{\Nsc}{\Leftrightarrow}
\newcommand{\qNsc}{\quad\Nsc\quad}
\newcommand{\fa}[1]{\forall#1\ }
\newcommand{\ex}[1]{\exists#1\ }
\newcommand{\fain}[2]{\forall{#1\in#2}\ }
\newcommand{\exin}[2]{\exists{#1\in#2}\ }
\newcommand{\fale}[2]{\forall{#1\leq#2}\ }
\newcommand{\exle}[2]{\exists{#1\leq#2}\ }

\newcommand{\lang}{L}
\newcommand{\Ltwo}{\lang_2}
\newcommand{\X}{\mathcal{X}}
\newcommand{\Y}{\mathcal{Y}}
\newcommand{\tuple}[1]{\langle#1\rangle}
\newcommand{\ind}{\mathnoun{I}}
\newcommand{\bd}{\mathnoun{B}}
\newcommand{\WKL}{\mathsf{WKL}}

\newcommand{\Def}{\mathnoun{Def}}
\newcommand{\len}{\mathnoun{len}}

\newcommand{\Dom}{\mathnoun{Dom}}
\newcommand{\M}{\mathfrak{M}}
\newcommand{\RCA}{\mathsf{RCA}_0}
\newcommand{\RT}{\mathsf{RT}^2_2}
\newcommand{\COH}{\mathsf{COH}}
\newcommand{\TT}{\mathsf{TT}^1}

\renewcommand{\exp}{\mathnoun{exp}}
\let\Tjoin\oplus

\newcommand{\defined}{\mathclose\downarrow}
\newcommand{\defm}[1]{\emph{#1}}

\theoremstyle{definition}

      \newtheorem{theorem}{Theorem}[section]
    \newtheorem{lemma}[theorem]{Lemma}
    \newtheorem{corollary}[theorem]{Corollary}

      \newtheorem{proposition}[theorem]{Proposition}

     \newtheorem{definition}[theorem]{Definition}
     \newtheorem{question}[theorem]{Quesition}
     \newtheorem*{rmk}{Remark}

\title{Definability over $\bd\Sigma^0_2$-models}

\author{Chi Tat Chong}
\address{Department of Mathematics, National University of Singapore\\
 Singapore 119076}

\author{Tin Lok Wong}
\address{Department of Mathematics, National University of Singapore\\
 Singapore 119076}

\thanks{Chong's research was partially supported by NUS grants C-146-000-042-001.}

\dedicatory{---To Qi Feng on his 70th birthday}

\begin{document}

\begin{abstract}
Let $\M=(M, \X)$ be a model of  $\mathsf{RCA}_0+\Sigma^0_2\text{-bounding}$ in which $\Sigma^0_2(A)$-induction fails for some  $A\in\X$. We show that (i) if $\M$ is a model of the combinatorial principle Ramsey's Theorem for Pairs, the Cohesive Set Theorem or the Tree Theorem, then there is a $\Delta^0_1(A)$-instance of the  principle  with no solution in $\M$ that is arithmetically definable  relative to $A$;
and  (ii)  any set of minimal Turing degree in $\M$  that is arithmetically definable relative to $A$ has Turing jump equivalent to $A'$.

\end{abstract}

\maketitle

 \section{Introduction}

 Definability is of central interest in the study of computation theory. Calibrating the complexity of a set constructed to solve a problem  often  brings insights  to  the  nature of the problem and its solution. This complexity may be measured in terms  of   definability in the language of  the model of computation.
 The  intuition is that existence or nonexistence of a definable solution
 offers a deeper understanding not only of the problem but also of the underlying computation model.
   The general question is, given a problem that is definable in the model, does it have a definable solution over the model?

 Classically there are countless  examples  illustrating this question over the standard model $
 \mathfrak N=(\mathbb N,  0,1, +, \times)$ of  arithmetic. We cite some prominent ones  which motivated  the study discussed in this paper.  These  relate to  Ramsey-type combinatorial principles and the Turing degree structure---more  precisely, the existence of definable solutions for instances of Ramsey's Theorem for Pairs  ($\RT$), the Cohesive Set Principle ($\COH$)  and the Tree Theorem Principle  ($\TT$)  (see \S \ref{Comb} for a definition of these notions),
  as well as the existence of a definable set of minimal Turing degree:

   \begin{enumerate}
 \item   Every recursive  two-coloring of pairs in $\mathbb N$  has a
 $\Pi^0_2$-definable homogeneous subset
   (Jockusch \cite[Theorem 4.2]{Jockusch});

 \item  Every recursive array $\{A_i: i\in\mathbb N\}$ of sets $A_i\subset\mathbb N$ has a cohesive set $A$ such that $A''\le_T \emptyset''$ (Jockusch and   Stephan \cite[Theorem 2.5]{JS});

 \item Every recursive coloring of the nodes in the full binary tree in finitely many colors has a recursive homogeneous subtree isomorphic to the full binary tree (immediate);

 \item There is a set  of minimal degree  $<_T\emptyset''$ (Spector \cite{Spector56});
\item There is a set of minimal degree $<_T \emptyset'$ (Sacks \cite{Sacks63}).
  \end{enumerate}

If $\M=(M, \X)\models\RCA+ \Sigma^0_2\text{-induction}$ ($\ind\Sigma^0_2$), then every instance $A\in\X$  of $\COH$ or $\RT$ has an $(M, A)$-definable  solution $G$
such that  $\M[G]\models \ind\Sigma^0_2$ by Ikari \cite[Theorem 4.1]{Ikari}.
And for $\TT$ there is always a recursive solution. In the case of  the minimal degree problem, the construction  in \cite{Sacks63} can be adapted to produce a $G<_T\emptyset'$ of minimal degree  and so $\M[G]\models\ind\Sigma^0_2$.

 In this paper  we consider  the question of definable solutions for examples (1)--(5)  above   in  the context of  models $\M=(M, \X)$  of $\RCA+\Sigma^0_2\text{-bounding}$ ($\bd\Sigma^0_2$)  in which $\ind\Sigma^0_2$ fails. This  is
  the weakest system extending $\RCA$, in terms of inductive strength,   in second-order arithmetic with a model  for each of $\RT$ and $\TT$ (cf.~Chong, Slaman and Yang \cite{CSY2017} and Patey and Yokoyama \cite{PY},  and Chong, Li, Wang and Yang \cite{CLWY} respectively), and a model for $\COH$
  which can be expanded  (by adding appropriate second-order elements)  to that of $\RT$  (Chong, Slaman and Yang \cite{CSY2012}).
     In the case of the minimal degree problem, the system is a natural analog in second-order arithmetic  of  $\Sigma_2$-inadmissible ordinals such as $\aleph^L_\omega$ or $\aleph^L_{\omega_1}$ where the existence of a minimal degree  has remained unknown since the problem was posed by Saacks  morfe than five decades ago.

 The construction   in \cite{CSY2017}, \cite{CSY2012} and \cite{CLWY} to obtain a  model of $\RT, \COH$ or  $\TT$ is carried out on a countable $\M$  by iteratively adding  sets  $G\subseteq M$ such that (i) $G$ is a solution of an instance of  the combinatorial principle concerned, and
     (ii)  $\Sigma^0_2$-bounding  is preserved upon adding  $G$.   While the countability of $\M$ ensures that  the construction will succeed, it does not exhibit  the  definability of $G$ over $\M$.  Indeed by Kossak \cite[Corollary
   3.3]{Ko} there is no countable model of this system  whose countability is witnessed definably over the model.     Hence  a definable solution for the  problems considered  can only be achieved through a  direct construction definable over the model.

   Following the preliminaries in Section~\ref{s:prelim}, we prove  in Section~\ref{s:Sigma03}  (Theorem~\ref{Sigma3})   that if $\M=(M, \X)\models\RCA+\bd\Sigma^0_2$ and $\ind\Sigma^0_2(A)$ fails for $A\in\X$, then every  $\Sigma^0_3(A)$-definable $G\subset M$ such that  $\M[G]\models\bd\Sigma^0_2$  is low relative to $A$, i.e.,~the Turing jump of $G$ is Turing equivalent to $A'$. In Section~\ref{s:RCA0*} we transition to the subsystem $\RCA^*+\neg\ind\Sigma^0_1$ (see \S\ref{s:prelim} for the definition)  and prove, for $n=0$ in Theorem \ref{4-2},  that if $(M, A)\models \RCA^*+\neg\ind\Sigma^0_1$ and if $B$ is arithmetically definable in $(M, A)$ such that $(M, A, B)\models  \bd\Sigma^0_1$, then $B$ is $\Delta^0_1(A)$, i.e.,~$B\le_T A$.
     As applications we prove in Section~\ref{s:app} that if $\M=(M, \X)\models \RCA+\bd\Sigma^0_2+\neg\ind\Sigma^0_2$   then over $\M$,  (i) there exist instances of $\RT, \COH$ and $\TT$ with no arithmetically definable solution relative to  the  instance, and (ii) there is an $A\in\X$ such that every set of minimal Turing degree arithmetically definable relative to $A$ is low relative to $A$, i.e.~has Turing jump equivalent to $A'$.
      We conclude the paper with a list of questions in Section~\ref{s:concl}.

 \section{Preliminaries}\label{s:prelim}

 \subsection{Notations, second-order arithmetic and recursion-theoretic notions}

  We work in  the language of second-order  arithmetic.  For $n\ge 0$, define a formula in the language to be  $\Sigma^0_n$  or  $\Pi^0_n$ as  usual (however, where appropriate we will use $\Sigma_n$ and $\Pi_n$ in   the first-order language of arithmetic, possibly with number parameters, to  emphasize the absence of second-order sets and quantifiers).   If $\M=(M, \X)$ is a structure in the language  of second-order arithmetic and $G\subseteq M$, then $G$ is $\Sigma^0_n(A)$,
   for some $A\in\X$,  if $G $ is $\Sigma^0_n$-definable over $\M$ with parameter  $A$ (and no other set parameters). Define ``$G$ is $\Pi^0_n(A)$''  similarly. We say that $G$ is
     $\Delta^0_n(A)$ for $A\in\X$  if $G$ is both $\Sigma^0_n(A)$ and $\Pi^0_n(A)$.

     Given $A\subseteq M$, let $(M, A)$ be the structure whose first-order part is $M$ and second-order part is $\{X: X\in\Delta^0_1(A)\}$.
  We write $G\in\Sigma^0_n(M, A)$ (or ``$G$ is $\Sigma^0_n(M, A)$'')  if it is  $\Sigma^0_n$-definable over $(M, A)$.\footnote{If $\M\models \RCA^*$ (defined below) then  $\Sigma^0_n(M,A)=\Sigma^0_n(A)$. This is the case for all models considered in this paper.}
    Define $G\in\Pi^0_n(M, A)$ similarly.
   If $\M=(M, \X)$,  then $G\in\Sigma^0_n(\M)$ if it is $\Sigma^0_n(M, A)$ for some $  A\in\X$.
  Define  $G\in \Pi^0_n(\M)$ and $G\in\Delta^0_n(\M)$ similarly. We say that $G$ is definable in $\M$ or in $(M, A)$ if it is $\Sigma^0_n(\M)$ or $\Sigma^0_n(M, A)$ for some $n$. Finally, if $G\subseteq M$, then $\M[G]=(M, \Y)$ where $\Y$ is the ideal generated by $\X$ and $G$, i.e.~$\Y=\{B: \exists A\in\X(B\in\Delta^0_1(A\oplus G))\}$.

  Let $P^-$ denote  formalization of  the Peano axioms minus the induction scheme as defined in Paris and Kirby \cite{PK}.
 For $n \ge 1$, $\ind\Sigma^0_n$ denotes the $\Sigma^0_n$-induction scheme and $\bd\Sigma^0_n$
   the $\Sigma^0_n$-bounding scheme.
     Over $P^- +\ind\Sigma_0$, $\bd\Sigma^0_{n+1}$ has proof-theoretic strength strictly between $\ind\Sigma^0_{n+1}$ and $\ind\Sigma^0_n$  \cite{PK}).
     Since $\bd\Sigma^0_n$ is equivalent to $I\Delta^0_n$  for  $n\ge 1$  by Slaman \cite{Slaman}, the chain
     $$
     \cdots\rightarrow \ind\Sigma^0_{n+1}\rightarrow \bd\Sigma^0_{n+1}\rightarrow \ind\Sigma^0_n\rightarrow\cdots
     $$
      is a natural ordering of arithmetical subsystems in terms of inductive strength.

      If $\M=(M, \X)$ is a structure  in   second-order arithmetic,  then $\M\models\exp$ asserts the totality of the  exponential function.
       It is immediate that $\exp$ follows from $\Sigma_1$-induction.

The domain of a function $f$ is denoted $\Dom(f)$.       A bounded subset of $M$ is \defm{$\M$-finite} if it is coded in $M$ using the exponential function.
A  \defm{(binary) string} $\sigma$  in $\M$ is an $\M$-finite  $\{0, 1\}$-valued function such that $\Dom(\sigma)=\text{an initial segment of } M$. The \defm{length} of $\sigma$, denoted  $\len(\sigma)$,  is the least  $x$ such that $\sigma(x)$ is not defined.

      There is a  well-developed theory of computation `over'' (first-order) models of $P^-+\ind\Sigma_n$  and $P^-+\bd\Sigma_n$ for $n\ge 1$  (with $\ind\Sigma_0+\exp$ as base theory when $n=1$),
             This theory carries over to
              models  of a sufficiently strong subsystem of second-order arithmetic. Recall that $\RCA$ is the subsystem of second-order arithmetic comprising $P^-$,  $\Sigma^0_1$-induction and   $\Delta^0_1$-comprehension.
   The system $\RCA^*$  is a weakening of $\RCA$,  in  that the  $\Sigma^0_1$-induction scheme is replaced by  $\Delta^0_1$-induction plus  $\exp$.

   Given $\M=(M, \X)$ a model of $\RCA^*$  and
   $A, B\subseteq M$, we say that  $A$ is {\it  recursive} (or \defm{computable}) in $B$, written $A\le_T  B$,  if $A$ is  $\Delta^0_1(B)$. Thus there is  a total $\Sigma^0_1(M, B)$-function  $\Phi$ whose output is (the characteristic  function of)  $A$. We denote this  by  $\Phi^B=A$. The set $A$ is {\it recursive}  (or \defm{computable}) if it is $\Delta_1$ over $\M$ (without set parameters).
              We say that  $A, B$ have the \defm{same Turing degree}, written $A\equiv_T B$, if $A\le_T B$ and $B\le_T A$.    Denote by $A<_T B$ if $A\le_T B$ and $B\not\le_T A$.
    The relation  $\le_T$ is transitive over models  that satisfy $\bd\Sigma^0_1$ and hence
    models of  $\RCA^*$.

    Suppose $\M= (M, \X)\models\RCA$. A notion that is particularly relevant  to our study in this paper is that of the \defm{Turing jump} $A'$ of an $A\in\X$.  It is defined to be the set $ \{e: \Phi^A_e(e) \defined\}$,  where $\Phi^A_e(e)\defined$ means $ \Phi^A_e(e)$  is defined  ($\Phi_e$ is the $e$th partial $\Sigma^0_1(M, A)$-function). While $(M, A')$ may not satisfy $\ind\Sigma^0_1$ and so  $A'$ may not be in $\X$, it is $\Sigma^{0}_1(M, A)$ and its Turing degree is well-defined in the sense that $G\subseteq  M$ is $\Delta^0_1(A')$ if and only if $G$ is $\Delta^0_2(A)$.\footnote{For
       $n\ge 1$  and
     $\M\models \bd\Sigma^0_{n+1}$,  the Turing degree of $A^{(n)}$ is well-defined, i.e.~$G\subseteq M$ is $\Delta^0_1(A^{(n)})$ if and only if $G$ is $\Delta^0_{n+1}(A)$ (Post's Theorem). Furthermore,
      $(M,  A^{(n)})\models \bd\Sigma^0_{1}$. }
      Furthermore,  since $(M, A)\models \ind\Sigma^0_1$, the set $A'$ is \defm{regular},  i.e.~$A'\restriction s$ is $\M$-finite for all $s\in M$.

         This observation points to  a natural connection between $\RCA+\bd\Sigma^0_2$ and $\RCA^*$, namely if $\M=(M, \X)$ is a model of $\RCA+\bd\Sigma^0_2$,
 then $\M'_A=(M, A')$   is a model of $\RCA^*$.
      (see Lemma \ref{BSigma2}(iv) below).
 The  final  subsystem of second-order arithmetic considered here is $\WKL^*_0$, which is $\RCA^*+
     \mathsf{WKL}$.  $\WKL$ is weak K\"onig's lemma, which states that every unbounded  $0$--$1$ (binary)  tree has an unbounded path. A well-known result due to Harrington states that every countable model of $\RCA$  can be extended to a model of $\RCA+\WKL$ with the same first-order universe. Similarly, a countable model of $\RCA^*$ can be extended to that of $\WKL^*_0=\RCA^*+ \WKL$ by adding only second-order elements~\cite[Theorem~4.6]{art:ss}. The argument in Simpson's book~\cite[Theorem IV 4.4]{Simpson} shows that $\WKL$ is equivalent over $\RCA^*$ to the $\Sigma^0_1$-separation principle, which states that every disjoint pair $(A, B)$ of $\Sigma^0_1$-sets is separated by  a set in the second-order universe---one which contains~$A$ and is disjoint from $B$. This fact will be used in the proof of  Theorem~\ref{4-1}. The next result of  Fiori-Carones, Ko\l odziejczyk, Wong and  Yokoyama \cite{CKWY} will be applied in the proof of Theorem \ref{4-2}:

     \begin{proposition}[{\cite[Theorem 2.1]{CKWY}}]\label{Iso}
If $(M, \X)$ and $(M,\Y)$ are countable models of $\RCA^*$ and $(M, \X\cap\Y)\models \neg\ind\Sigma^0_1$, then $ (M, \X)$ and $(M, \Y)$ are isomorphic.  Furthermore, the isomorphism can be required to fix any prescribed finite tuple of first- or second-order parameters from $(M,\X\cap\Y)$.

   \end{proposition}

    \vskip.15in

             \subsection{$\bd\Sigma^0_2$-models and  $\bd\Sigma^0_2$-sets}
 \begin{definition}

  Suppose   $\M=(M,\X)\models\RCA+\bd\Sigma^0_2$.
  \begin{enumerate}
  \item [(i)]  $\M$  is a \defm{$\bd\Sigma^0_2$-model} if  $\M\models \neg \ind\Sigma^0_2$ (if $\M\models\RCA^*+\neg\ind\Sigma^0_1$ then $\M$ is a $\bd\Sigma^0_1$-model).

  \item [(ii)] $G\subset M$ is a \defm{$\bd\Sigma^0_2$-set} over $\M$  if $\M[G]\models \RCA+\bd\Sigma^0_2$ ($G$ is a $\bd\Sigma^0_1$-set if $\M[G]\models\RCA^*$).
 \end{enumerate}

 \end{definition}

   The following   are basic facts about  $\bd\Sigma^0_2$-models $\M=(M, \X)$:

 \begin{itemize}
 \item $\M$ has a $\Sigma^0_2$-cut: There is  an $A\in\X$ and a bounded, downward-closed  $\Sigma^0_2(A)$-definable   $I\subset M$ which  is closed under the successor operation.
 \item  There is an increasing  $\Sigma^0_2(A)$-definable  function $g$ with $\Dom(g)=I$ and range cofinal in~$M$.
 \item  There is an approximation $g'\in\Delta^0_1(A)$  of $g$ such that
 \begin{enumerate}
 \item [(i)]  $\Dom(g')=M\times [0,a]$ for some $a$ which is an upper bound of $I$;
 \item [(ii)] For all $x\le y\le a$ and $s\le t$, $g'(s,x)\le g'(t,y)$;
 \item [(iii)]  $\text{lim}_s \ g'(s,x) $ exists   if and only if  $x\in I$, and in which case it is equal to $g(x)$;
 \item [(iv)]  The graph of $g$ is $\Delta^0_1(A')$;
 \end{enumerate}
 \item (Limit Lemma)  Every   $\Sigma^0_2(A)$-function $h$ that is $\{0,1\}$-valued has a   primitive recursive (relative to $A$)  approximation  $h': M\times M\rightarrow \{0,1\}$ such that
 $$
 \text{lim}_s\ h'(s, x)\text{ exists } \Leftrightarrow x\in\Dom(h),
 $$
 in which case $\text{lim}_s\ h'(s, x)= h(x)$.
  \end{itemize}

 \vskip.1in
 Given a $\bd\Sigma^0_2$-model $\M=(M, \X)$ and a bounded set $X\subset M$, we say that $Y\subseteq X$ is \emph{$\Delta^{0}_2(\M)$ on $X$} if both $Y$ and $X\setminus Y$ are $\Sigma^{0}_2(\M)$. We say that  $Y$ is {\it coded on} $X$ if there is an $\M$-finite set $\hat Y$ such that $Y=X\cap\hat Y$. The  next proposition is Proposition~4 of  Chong and Mourad~\cite{CM}  cast in the setting  of second-order arithmetic.

  \begin{proposition}\label{CM}
 Let $\M=(M, \X)$ be a $\bd\Sigma^0_2$-model. If $X$  is a bounded subset of $M$ and $Y\subseteq X$ is $\Delta^{0}_2(\M)$ on $X$,  then    $Y$ is coded on $X$.
 \end{proposition}

 \begin{lemma}\label{BSigma2}
 Let $\M=(M, \X)$ be a $\bd\Sigma^0_2$-model with a $\Sigma^0_2(A)$-cut $I$ for some $A\in\X$, and let $G$ be a  $\bd\Sigma^0_2$-set over $\M$. Then
 \begin{enumerate}\renewcommand{\theenumi}{\roman{enumi}}
\item $G$ is regular;
 \item $G\oplus A\not\ge_T   A'$ and $G\oplus A\not\ge_T I$;
 \item $G\oplus A$ is $A\text{-}\mathsf{GL}_1 $  (generalized low relative to $A$), i.e.~$(G\oplus A)'\equiv_T G\oplus A'$;\label{pt:BSigma2:GL1}
 \item For each $A\in\X$, $\M'_A=(M, A')$ is a $\bd\Sigma^0_1$-model and $(G\oplus A)'$  is a $\bd\Sigma^0_1$-set over  $\M'_A$.
 \end{enumerate}
 \end{lemma}

 \begin{proof}
 (i).  Given $s\in M$, the set $\{s: \text{$s'\le s$ and  $s'\in G$}\}$ is $\Sigma^0_1(G)$ and hence $\M$-finite by $\ind\Sigma^0_1(G)$. Hence $G$ is regular.

 (ii). Since $G$ is a $\bd\Sigma^0_2$-set, we have $\bd\Sigma^0_2$ to hold relative to $G\oplus A$ as well. Now suppose $G\oplus A\ge_T  A'$.  Let $g\colon I\rightarrow M$ be $\Sigma^0_2(A)$, strictly increasing and cofinal.  Then $g\in \Delta^0_1(A')$ by basic fact (iv) above and hence  is $\Delta^0_1((G\oplus A)')$.  Thus  there is a partial $\Sigma^0_1$-function $\Phi$ with a free set variable such that  $\Dom(\Phi^{(G\oplus A)'})=\Dom(g)=I$.  This contradicts $\ind\Sigma^0_1((G\oplus A)')$.

   (iii).  Let $\Phi^{G\oplus  A}$ define $(G\oplus A)'$. By $\ind\Sigma^0_1$ relative to $G\oplus A$, for every   $i\in I$, $\Dom(\Phi^{G\oplus A})\restriction g(i)$ is $\M$-finite. By  $\bd\Sigma^0_1$ for $G\oplus A$, there is a $j\in I$ such that $\Phi^{G\oplus A} \restriction g(i)=\Phi^{(G\oplus A)\restriction g(j)}\restriction g(i)$.  The set
 $$
 Y= \{(i,j): j\text{ is the least } j' \text{  such that }\Phi^{G\oplus A}\restriction g(i)=\Phi^{(G\oplus A)\restriction g(j')}\restriction g(i)\}
 $$
 is $\Delta^0_2(M, G\oplus A)$ on $I\times I$ and hence by Proposition \ref{CM} is coded on $I\times I$ by an $\M$-finite set $\hat Y$.
 This  code provides an algorithm  to compute $(G\oplus A)'$ from $(G\oplus A)\oplus A'$: Given $a\in M$, use $A'$ to find the least $i\in I$ such that $a\le g(i)$. Let $j$ be chosen so that $(i,j)\in\hat Y$. Then  for $x\le g(i)$, $x\in (G\oplus A)'$ if and only if $\Phi^{(G\oplus A)\restriction g(j)}(x)=1$.
Hence $(G\oplus A)'\le_T  (G\oplus A)\oplus A'\equiv_T G\oplus A'$.  The proof that $G\oplus A'\le_T (G\oplus A)'$ is immediate.

 (iv). Immediate since $A$ is a $\bd\Sigma^0_2$-set.
 \end{proof}

\subsection{Combinatorial principles}\label{Comb}
We consider three combinatorial principles which have been studied in reverse mathematics: the cohesive set principle $\mathsf{COH}$, Ramsey's theorem for pairs $\RT$, and the tree theorem $\mathsf{TT}^1$. For the discussion that follows, we recall the definitions of these principles.

\begin{definition}
\begin{enumerate}\label{RamDef}
\item [(i)] $\COH$: Given an array $A_0, A_1,\dots, A_e, \dots$  of sets, there is an infinite $C$ such that for each $e$, either $C\setminus A_e$ is finite or  $C\cap A_e$ is finite. $C$ is  said to be \defm{cohesive} for the array.
\item [(ii)] $\RT$: Every two-coloring of pairs of numbers has an infinite homogeneous set $H$, i.e.~all pairs of numbers in $H$ have the same color.
\item [(iii)]  $\TT$: Every finite coloring of the full $0$--$1$-(binary) tree has an isomorphic homogeneous subtree, i.e.~a tree all of whose nodes have the same color.
\end{enumerate}
\end{definition}

In  the absence of $\Sigma^0_2$-induction,  there is   a completely different  picture compared with that shown in \cite{Ikari}.
The next lemma rules out a solution below $\emptyset'$.

\begin{lemma}\label{Below 0'}
Let $\M=	(M, \X)$ be a $\bd\Sigma^0_2$-model and let $A\in \X$.
For each  of $\COH, \RT$ and $\TT$  there is a  $\Delta^0_1(A)$-instance  for which    no   solution  $G$ satisfying $\M[G]\models \bd\Sigma^0_2$ is $\Delta^0_2(A)$-definable.
\end{lemma}

\begin{proof}
For $\RT$ and $\TT$, these follow respectively  from  Chong, Slaman and Yang \cite[Proposition 2.4]{CSY2014} and Chong, Li, Wang and Yang \cite[Theorem 2.2]{CLWY}.   We  prove  the lemma for $\COH$.

To simplify notations, suppose  $A=\emptyset$ (the general case is obtained by relativizing to $A$).  We
give   a $\Delta_1$-instance of $\COH$ with no $\Delta_2$-definable solution $G$ that is a    $\bd\Sigma^0_2$-set. Following  Jockusch and Stephan~\cite[Theorem~2.1]{JS}\footnote{We thank Frank
Stephan for alerting us to \cite{JS}.},  let $\{A_0, A_1,\dots\}$ be an array consisting of all  primitive recursive sets and suppose $G$ is  a $\bd\Sigma^0_2$-set that is   $\Delta_2$-definable  and cohesive for the array. We claim  that $G'>_T \emptyset'$. This will yield a contradiction since by
Lemma \ref{BSigma2}(i) a $\Delta_2$-definable  set that is $\bd\Sigma^0_2$ has to be below $\emptyset'$ in Turing degree, and hence by Lemma \ref{BSigma2}(ii) we must have $G'\equiv_T\emptyset'$.
Let  $h$ be a partial $\{0, 1\}$-valued $\Sigma_2$-function.

\vskip.15in
{\bf  Claim 1.} There is a total $\{0,1\}$-valued function $\hat h$ which is $\Delta_1(G')$ and extends $h$.

By Basic Fact (v) let $h'\colon M\times M\rightarrow \{0, 1\}$ be a primitive recursive approximation of $h$. Hence
$$
\text{lim}_s\ h'(s,x) \text{ exists   }  \Leftrightarrow\ x\in\Dom(h),
$$
and $\text{lim}_s\ h'(s,x)$  exists implies that the limit is $h(x)$.
Then there is a $\Delta_1$-definable  function $x\mapsto\hat x$ such that for all $x$, $A_{\hat x}=\{s: h'(s,x)=1\}$.
Now since $G$ is cohesive for the array of primitive recursive sets, $\text{lim}_{s\in G}\ h'(s, x)$ exists for each $x$,  as either $\{x: G(x)=1\}\setminus A_{\hat x}$ or $\{x: G(x)=1\}\cap A_{\hat x}$ is $\M$-finite.

Since $\M[G]\models \ind\Sigma^0_1$, for each $x$ there is an $s(x)$ with the property that
 \begin{align*}
 s(x)=\text{ the least } s \text{ such that }&\forall s'\ge s\ (s'\in G\rightarrow s'\in   A_{\hat x})\text { or }\\                    &\forall s'\ge s\ (s'\in G\rightarrow s'\notin A_{\hat x}).
                    \end{align*}

 It follows that there is a total $\{0,1\}$-valued $\Delta_1(G')$-function $\hat h \supset h$ such that for all $x$, $\hat h(x)=1$ if and only if  $G\subseteq A_{\hat x}$ above $s(x)$ and $\hat h(x)=0$ if and only if  $G\cap A_{\hat x}=\emptyset$ above $s(x)$.

  The uniformity of transiting from $h$ to $\hat h$ shows  that  there is a $\Sigma^0_1$-function  $\Phi$ with a free set variable   such that  for each $e$, $\Phi^{G'}(e)=\hat e$ is an index of  the total $\{0, 1\}$-function $\hat{h_e}$ which extends the  partial $\{0,1\}$-valued  $\Sigma_2$-function $h_e$.

  \vskip.15in
  {\bf Claim 2.} $G'>_T\emptyset'$.

  Suppose $G'\equiv_T\emptyset'$. Then by $\Phi$ above every partial $\Sigma_2$-function is uniformly extendible to a total $\Delta^0_1(\emptyset')$-function. Define
  \begin{equation*}
  \theta (e)=
  \begin{cases}
  0,  &\text{if }  \hat h_{ e}(\hat e)=1;\\
1,&\text{otherwise}.
\end{cases}
\end{equation*}
Then $\theta$ is $\Delta^0_1(\emptyset')$ and hence $\Delta^0_2$. Let $\theta=h_{e_0}$ for some $e_0$. Since $\theta$ is total, we have $h_{e_0}=\hat h_{e_0}$.  But then $\theta(e_0)=1$ if and only if $\theta(e_0)\ne 1$, which is a contradiction. This completes the proof of the lemma.
\end{proof}

 \section{ $\Sigma^0_3$-definable $\bd\Sigma^0_2$-sets}\label{s:Sigma03}

 Let $\M$ be a $\bd\Sigma^0_2$-model.
  In this section, we analyze ``internally'', i.e.~without reference to the cardinality of $\M$  (as opposed to ``externally'' in the next section),  the question of  existence of definable $\bd\Sigma^0_2$-sets over $\M$.

\begin{theorem}\label{Sigma3}
Let $\M=(M, \X)$ be a $\bd\Sigma^0_2$-model, and $A\in\X$ for which $\ind\Sigma^0_2$ fails.
If $G\subset M$  is $\Sigma^0_3(A)$
 and $\M[G]\models \bd\Sigma^0_2$, then   $G'\le_T  A'$.
\end{theorem}

\begin{proof}
Fix a $\Sigma^0_2(A)$-cut $I$ of~$M$ and
 a $\Sigma^0_2(A)$-definable increasing and cofinal function
  $g\colon I\to M$.
Let
 \begin{equation*}
  G=\{x\in M : \M\models\ex u\fa v\ex w\varphi(x,u,v,w,A)\},
 \end{equation*} where $\varphi\in\Delta^0_0$.
Consider the set
 \begin{multline*}
  C=\{(i,j,k)\in I^3 :
   \M\models\fale x{g(k)}\\\bigl(
    x\notin G\then
    \fale u{g(i)} \exle v{g(j)} \fa w \neg\varphi(x,u,v,w,A)
   \bigr)
  \}.
 \end{multline*}
Then $C$ is coded on~$I^3$ by Lemma~\ref{CM}.
Find an $\M$-finite set $\hat C$ such that $C=I^3\cap\hat C$.
We claim that
 \begin{multline*}
  G=\bigl\{x\in M : \M\models\exin{i,k}I \exle u{g(i)} \fain jI\\
   \bigl(
    (i,j,k)\in\hat C
    \then g(k)\geq x\wedge\fale v{g(j)} \ex w \varphi(x,u,v,w,A)
   \bigr)
  \bigr\}.
 \end{multline*}
An application of $\bd\Sigma^0_1(A)$ will then tell us
 $G$ is~$\Sigma^0_2(A)$.
By repeating all these to
  the resulting $\Pi^0_2(A)$-definition of $M\setm G$,
 we will be able to deduce that $G$ is~$\Delta^0_2(A)$,
  or $G\le_T A'$.
From this, we will get
 $G'\le_T(G\Tjoin A)'\equiv_T G\Tjoin A'\le_T A'$
  by Lemma~\ref{BSigma2}\ptref{pt:BSigma2:GL1}, as required.

For the left-to-right direction of the claim, take any $x_0\in G$.
Then
 \begin{equation*}
  \M\models\ex u\fa v\ex w\varphi(x_0,u,v,w,A)
 \end{equation*} by the definition of~$G$.
Using the cofinality of~$g$ to find $i,k\in I$ and $u_0\leq g(i)$
 such that
 \begin{math}
  \M\models\fa v\ex w\varphi(x_0,u_0,v,w,A)
 \end{math} and $g(k)\geq x_0$.
Pick any $j\in I$ satisfying $(i,j,k)\in\hat C$.
Note that we have
 \begin{math}
  \M\models\fale v{g(j)} \ex w \varphi(x_0,u_0,v,w,A)
 \end{math},
 because this is true even without the bound~$g(j)$.
Hence these $i,k$ satisfy the requirements.

For the right-to-left direction,
 take any $x_0\in M\setm G$.
Then
 \begin{equation*}
  \M\models\fa u\ex v\fa w\neg\varphi(x_0,u,v,w,A)
 \end{equation*} by the definition of~$G$.
Pick any $i,k\in I$ and $u_0\leq g(i)$
 such that $g(k)\geq x_0$.
By the definition of~$G$ again, we know
 \begin{equation*}
  \M\models\fale x{g(k)}
   \bigl(
    x\notin G\then
    \fale u{g(i)} \ex v \fa w \neg\varphi(x,u,v,w,A)
   \bigr).
 \end{equation*}
Replacing $G$ with the $\M$-finite set $G\restriction g(k)$ here,
 we get from $\bd\Pi^0_1(A)$ and the cofinality of~$g$
  some $j\in I$ such that
  \begin{equation*}
   \M\models\fale x{g(k)}
    \bigl(
     x\notin G\then
     \fale u{g(i)} \exle v{g(j)} \fa w \neg\varphi(x,u,v,w,A)
    \bigr).
  \end{equation*}
Then $(i,j,k)\in C\subseteq\hat C$ and
 \begin{math}
  \M\models \exle v{g(j)} \fa w \neg\varphi(x_0,u_0,v,w,A)
 \end{math}.
\end{proof}

\begin{rmk}
It follows immediately from Theorem \ref{Sigma3} that every $\Pi^0_3(A)$-definable $\bd\Sigma^0_2$-set is strictly below $A'$ in
Turing degree. Furthermore, the argument  shows that if  $G$ is regular and not necessarily a $\bd\Sigma^0_2$-set, then $G$ is $\Delta^0_1(A')$.
\end{rmk}

 \section{Over $\mathsf{RCA}^*_0$}\label{s:RCA0*}

 In this section, we generalize Theorem \ref{Sigma3}  on $\Sigma^0_n(A)$-definable $\bd\Sigma^0_2$-sets  to all $n\ge 2$.  We achieve  this by first  constructing a model of  $\WKL^*_0$  with restricted collection of definable second-order elements:

\begin{theorem}\label{4-1}
Let $\M=(M,\X)$ be a countable model of $\RCA^*$, and let
 $A\in\X$.
Then there exists a $\Y$ such that  $(M,\Y)\models\WKL_0^*$ and for any $Y\in\Y$,
$$
Y\text{ is second-order definable over  } (M, \X)\Leftrightarrow Y\in\Delta^0_1(A).
$$
 \end{theorem}

\begin{proof}
We use the approach in the proof of  \cite[Lemma~3.2]{CKWY}.
Roughly speaking, take elements of~$M$ to Cohen-approximate a subset of $M$
 whose columns form our~$\Y$.
The set of all such elements is
 an $M$-cofinal tree~$T_A$ that is $\Delta^0_1$-definable over $(M,A)$,
  \emph{any} path through which consists
   of columns that form a model of~$\WKL_0^*$ containing~$A$.
   This is guaranteed by having designated columns to witness
 the containment of~$A$ (see~(\ref{cl:T:A}) below),
 closure under Turing join (see~(\ref{cl:T:join1}) and~(\ref{cl:T:join2}) below), and
 $\Sigma^0_1$~separation (to  satisfy  $ \WKL^*_0$; see~(\ref{cl:T:sep}) below).
We then choose a path through~$T_A$ carefully so that
 the resulting model of~$\WKL_0^*$ satisfies the required non-definability condition.

More precisely,
 consider the tree $T_A$ of all $\M$-finite binary strings   $\sigma\in2^{<{M}}$
 satisfying all of the following conditions.
\begin{enumerate}
\item\label{cl:T:A}
For all $\tuple{0,x}\in\Dom(\sigma)$,
 \begin{equation*}
  \sigma(0,x)=1\qNsc x\in A.
 \end{equation*}
\item\label{cl:T:join1}
For all $\tuple{2\tuple{i,j}+1,2x}\in\Dom(\sigma)$,
 \begin{equation*}
  \sigma(2\tuple{i,j}+1,2x)=1\qNsc\sigma(i,x)=1.
 \end{equation*}
\item\label{cl:T:join2}
For all $\tuple{2\tuple{i,j}+1,2x+1}\in\Dom(\sigma)$,
 \begin{equation*}
  \sigma(2\tuple{i,j}+1,2x+1)=1\qNsc\sigma(j,x)=1.
 \end{equation*}
\item\label{cl:T:sep}
For all $\tuple{2\tuple{d,e,k}+2,x}\in\Dom(\sigma)$,
 \begin{itemize}
 \item  If
  \begin{math}
   (M\restriction\len(\sigma), +,\times) \models  \Phi_d^{\{y : \sigma(k,y)=1\}}(x)\defined
   \wedge \neg\Phi_e^{\{y : \sigma(k,y)=1\}}(x)\defined
  \end{math},
  then $\sigma(\tuple{2\tuple{d,e,k}+2,x})=1$.

 \item If
  \begin{math}
   (M\restriction\len(\sigma), +,\times)    \models  \neg\Phi_d^{\{y : \sigma(k,y)=1\}}(x)\defined
   \wedge  \Phi_e^{\{y : \sigma(k,y)=1\}}(x)\defined
  \end{math},
  then $\sigma(\tuple{2\tuple{d,e,k}+2,x})=0$.
 \end{itemize}
 \end{enumerate}
This is a $\Delta^0_1(A)$~definition of~$T_A$. It is in fact a $\Pi^0_1(A)$-class, in the sense that if  $\sigma\notin T_A$ then no extension of $\sigma$ will be in $T_A$. This property allows one to conclude, using $\bd\Sigma^0_1(A)$,  that if $K$  is an $\M$-finite collection of pairwise incompatible strings in $T_A$  in which every  member is extended by  at most  boundedly many strings in $T_A$, then there is a uniform bound on the lengths  of  strings extending some  member of $K$.

To show that~$T_A$ is $M$-cofinal,
 given $b\in M$, one can use induction on~$w<b$
  to prove the existence of $\sigma\in T_A$ of length~$w$.
  Since all quantifiers involved can be bounded by $2^b$ here,
 this argument can be achieved using $\Delta_0(A)$~induction.

Next, we build an increasing $M$-cofinal sequence $(\sigma_n : n\in\omega)$
  of elements of~$T_A$
 with the inductive condition that
  there are $M$-unboundedly many elements of~$T_A$ above~$\sigma_n$
   at each step $n\in\omega$.
Use the countability of~$\M$ to obtain a strictly increasing
 cofinal sequence $(b_n: n\in\omega)$ in~$M$.
 Let $\{S_\ell: \ell\in\omega\}$ be a complete list of all subsets of $M$ definable (in the language of second-order arithmetic with parameters) over $(M, \X)$.
 Define $\sigma_0$ to be the empty string.
Suppose $i,\ell\in\omega$
 and $\sigma_{\tuple{i,\ell}}$ satisfying the inductive condition
  is defined.
We seek $\sigma_{\tuple{i,\ell}+1}$
  of length greater than~$b_{\tuple{i,\ell}}$ to force
 the $i$th column of the set with characteristic function $\bigcup_{n\in\omega}\sigma_n$
  to be different from~$S_\ell$.
Use $\bd\Sigma_1(A)$ to extend $\sigma_{\tuple{i,\ell}}$
 to $\sigma_{\tuple{i,\ell}}'$ of length greater than~$b_{\tuple{i,\ell}}$
  above which there are $M$-unboundedly many elements of~$T_A$.
Now consider the following statement:
\begin{equation}\tag{\textasteriskcentered}\label{diag}
 \begin{minipage}[t]{.8\textwidth}\raggedright
 There exists $\tau\succeq\sigma_{\tuple{i,\ell}}'$
  with $M$-unboundedly many elements of~$T_A$ above it and
       $\tuple{i,x}<\len(\tau)$ in~$M$
  such that
  \begin{align*}
   \text{either}\quad&\text{$\tau(i,x)=1$ and $x\not\in S_\ell$},\\
   \text{or}\quad&\text{$\tau(i,x)=0$ and $x\in S_\ell$},
  \end{align*}
 \end{minipage}
\end{equation}

  If \eqref{diag} holds, then setting $\sigma_{\tuple{i,\ell}+1}$ to be any such~$\tau$ would do the job.
So suppose not,
 i.e.,~for every $\tau\succeq\sigma_{\tuple{i,\ell}}'$,
 if some $\tuple{i,x}<\len(\tau)$ makes
  \begin{equation*}
   \neg(\tau(i,x)=1\qNsc x\in S_\ell),
  \end{equation*}
  then there must be only $M$-boundedly many elements of~$T_A$ above~$\tau$.

We claim that
 \begin{align}\addtocounter{equation}{4}
  S_\ell\label{Sl}
  &=\bigl\{x\in M :  (M, +,\times)\models\ex b \fain\tau{T_A}\\ &\hspace{10em}\bigl(
    \tau\supseteq\sigma_{\tuple{i,\ell}}'\wedge
    \len(\tau)=b
    \then\tau(i,x)=1
   \bigr)
  \bigr\},\nonumber\\
\intertext{and}
  M\setm S_\ell\label{Slc}
  &=\bigl\{x\in M :  (M, +,\times)\models\ex b \fain\tau{T_A}\\ &\hspace{10em}\bigl(
    \tau\supseteq\sigma_{\tuple{i,\ell}}'\wedge
    \len(\tau)=b
    \then\tau(i,x)=0
   \bigr)
  \bigr\}.\nonumber
 \end{align}
This will imply $S_\ell\in\Delta^0_1(M,A)$, and
 thus we can set $\sigma_{\tuple{i,\ell}+1}$
  to be $\sigma_{\tuple{i,\ell}}'$.

Proving the right-to-left direction of the claim is straightforward.
We will show  the other  direction.
We content ourselves with a verification of~\eqref{Sl};
 a verification of~\eqref{Slc}  is similar.
Take $x\in S_\ell$.
Consider $n=\max\{\tuple{i,x}+1,\len(\sigma_{\tuple{i,\ell}}')\}$.
Let
 \begin{equation*}
  s_{n,0}=\{\tau\in T_A :
   \len(\tau)=n\wedge\tau\supseteq\sigma_{\tuple{i,\ell}}'\wedge\tau(i,x)=0
  \}.
 \end{equation*}
This set is coded in $\M$ by $\ind\Delta_0(A)+\exp$.
By our supposition that \eqref{diag} does not hold,
 for every $\tau\in s_{n,0}$, there exists $b_\tau\in M$
  such that no element of~$T_A$ of length~$b_\tau$ extends~$\tau$.
An application of $\bd\Sigma_1(A)$ then gives $b\in M$
 such that, for every $\tau\in s_{n,0}$,
  no element of $T_A$ of length~$b$ extends~$\tau$.
This $b$ satisfies the requirement.

Let $\Y$ be the columns of $\bigcup_{n\in\omega}\sigma_n$. It follows that $(M, \Y)$ satisfies the conclusion of the theorem.
\end{proof}

\begin{theorem}\label{4-2}
Fix $n\in\omega$.
Let $A,B\subseteq M$
 such that $(M,A,B)\models\bd\Sigma^0_{n+1}+\exp+\neg\ind\Sigma_{n+1}(A)$.
If $B$ is arithmetically definable in~$(M,A)$,
 then $B$ is $\Delta^0_{n+1}$-definable in $(M,A)$.
\end{theorem}

\begin{proof}
In view of L\"owenheim--Skolem, we assume $\M$ is countable.
Let $S_{n}^{A\Tjoin B}$ be
 the complete $\Sigma_{n}$-definable set
  relative to the Turing join of $A$ and~$B$.
   Then $(M,S_{n}^{A\Tjoin B})\models\bd\Sigma^0_1+\exp+\neg\ind\Sigma^0_1$.
By  Simpson and Smith~\cite[Theorem~4.6]{art:ss}
 one can  extend $(M,S_{n}^{A\Tjoin B})$ to a countable model $(M,\X)\models\WKL_0^*$.
Notice that $S^A_{n}$, the complete $\Sigma_{n}$-set relative to $A$, is a member of $\X$ and
 $(M, \X)\models\neg\ind\Sigma_1(S_{n}^A)$ (for
  $n\ge  2$, $S^{A}_n $ is Turing equivalent to $A^{(n)}$).
  Apply Theorem \ref{4-1}  to obtain another countable model
  $(M,\Y)\models\WKL_0^*$
 such that every second-order definable member of $\Y$ is $\Delta^0_1(M,S_{n}^A)$.
Proposition \ref{Iso}  then gives
 an isomorphism $\pi\colon(M,\X)\to(M,\Y)$
  fixing $S_{n}^A$ as a second-order object.

Now suppose $B$ is arithmetically definable in~$(M,A)$.
Then $B$ is arithmetically definable in~$(M,S_{n}^A)$
 because $A$ is arithmetically definable in~$(M,S_{n}^A)$.
So $\pi(B)$ is arithmetically definable in $(\pi(M),\pi(S_{n}^A))=(M,S_{n}^A)$.
Since $\pi(B)\in\Y$, and writing $\Ltwo\hyp\Def(M,\X)$ for the collection of sets second-order definable over $(M, \X)$, one has
 \begin{equation*}
  \pi(B)\in\Y\cap(\Ltwo\hyp\Def(M,\X))=\Delta^0_1(M,S_{n}^A)
  =\Delta^0_1(\pi(M),\pi(S_{n}^A)).
 \end{equation*}
Applying $\pi^{-1}$ to both sides, we deduce that
 $B\in\Delta^0_1(M,S_{n}^A)=\Delta^0_{n}(M,A)$.
\end{proof}

\begin{rmk}
As one can see from the proof above,
 we can weaken the condition that
  $(M,A,B)\models\bd\Sigma^0_{n+1}+\exp+\neg\ind\Sigma_{n+1}(A)$
 to $(M,A)\models\bd\Sigma^0_{n+1}+\exp+\neg\ind\Sigma^0_{n+1}$
  plus $(M,S_{n}^A,B)\models\bd\Sigma^0_1$. We will use this fact in the next section.
\end{rmk}

\section{Applications to definable solutions}\label{s:app}
We now apply results in Sections~\ref{s:Sigma03} and~\ref{s:RCA0*}  to characterize   definable solutions of combinatorial principles and sets of minimal degree  over $\bd\Sigma^0_2$-models.

\begin{corollary}\label{main1}
Let $\mathsf P$ be either $\COH, \RT$ or $\TT$ and suppose that  $\M\models \RCA+\bd\Sigma^0_2+\neg\ind\Sigma^0_2$. If $\M\models \mathsf P$ then there is an $A\in\X$ and a $\Delta^0_1(A)$-instance  of $\mathsf P$ with no arithmetically definable solution relative to $A$   in  $\M$. In other words, there is no solution of $C$ arithmetically definable over $(M, A)$.

\end{corollary}

\begin{proof}

Again by the L\"owenheim--Skolem theorem, we may assume that $\M=(M, \X)$ is countable.  Let $A\in\X$ be such that $\M\models\neg\ind\Sigma^0_2(A)$.
By Lemma~\ref{Below 0'} there is a $\Delta^0_1(A)$-instance $C$  of $\mathsf P$ with no $\Delta^0_1(A')$-definable solution.
Assume that $G\in \X$  is a   $\Sigma^0_n(A)$-definable solution of $C$  for some $n$. Then $n >3$  by Theorem~\ref{Sigma3}.
Then $(M, A)\models\RCA^*$ by Lemma \ref{BSigma2}(iv)
and $G$ is $\Sigma_{n-1}$-definable over $(M, A')$. Since $(M, A')\models \bd\Sigma^0_1+\neg\ind\Sigma^0_1$, by Theorem \ref{4-2} and the remark following it,  we must have $n=2$ so that $G\in\Delta^0_1(A')$. But  no such solution exists for  $C$.
 This contradiction proves the   corollary.
\end{proof}

A set  $G>_T \emptyset$ is of \defm{minimal Turing degree} if every $A<_T G$ is computable.
The constructions in  \cite{Sacks63}  over the standard model may be adapted to obtain, over any $\M\models\RCA+\ind\Sigma^0_2$,
 a set $G<_T\emptyset'$ of minimal degree
 such  that $\M[G]\models \ind\Sigma^0_2$.
By contrast, it is not known if there is an arithmetically definable set of minimal degree  for  a $\bd\Sigma^0_2$-model, although for a countable $\M$, a non-definable solution can be constructed  using existing technology.
The next corollary   addresses the question of   existence of a definable solution.

\begin{corollary}
Let $\M=(M, \X)$ be a $\bd\Sigma^0_2$-model  and let $A\in\X$ satisfy $\M\models\neg\ind\Sigma^0_2(A)$.  If $G\in\X$ is a set of minimal degree and arithmetically definable  relative to $A$, then $G'\equiv_T A'$. In particular,  if $A=\emptyset$ then $G$ is low, i.e.~$G'\equiv_T\emptyset'$.

\end{corollary}

\begin{proof}
The same argument as in Corollary \ref{main1} shows that $G\in\Delta^0_1(A')$. Since $A'$ is not of minimal degree, we must have $G<_T A'$.   Furthermore, since $G$ is a $\bd\Sigma^0_2$-set,  Lemma \ref{BSigma2}(iii) implies that  $G'\equiv_T A'$.
\end{proof}

In view of \cite[Theorem 4.1]{Ikari} and the above\footnote{We thank Keita Yokoyama for alerting us to~\cite{Ikari}.},  one arrives at the following characterization of definable solutions:

\begin{corollary}
Let $\M=(M, \X)\models\RCA+\bd\Sigma^0_2$. The following are equivalent:
\begin{enumerate}
\item [(i)] $\M\models \ind\Sigma^0_2$;
\item [(ii)]  For every $A\in\X$, every $\Delta^0_1(A)$-instance of $\COH,  \RT$ or $\TT$ has an arithmetically definable solution relative to $A$ that satisfies  $\bd\Sigma^0_2$;

\item [(iii)] For every $A\in\X$, there is a set of minimal degree $\not\le_T A'$, arithmetically definable relative to $A$,  and  satisfies  $\bd\Sigma^0_2$.
\end{enumerate}
\end{corollary}

\section{Concluding remarks and questions}\label{s:concl}

We end this paper with some remarks and a list of questions.

\begin{question}
Let $(M,\X)\models\ind\Sigma_1+\bd\Sigma^0_1+\neg\ind\Sigma^0_1$.
Partially order, by Turing reduction,
 those $A\in\X$ for which $(M, A)\models\neg\ind\Sigma_1$.
Can this poset have a minimal element?
\end{question}

\begin{question}
Are there $(M,A),(M,B)\models\bd\Sigma^0_1+\exp+\neg\ind\Sigma_1$
 such that $(M,A,B)\not\models\bd\Sigma^0_1$?
\end{question}

By  Theorem \ref{4-2},
 one way to get a positive answer to this question is to
  find $A,B\subseteq M\models\bd\Sigma_1+\exp+\neg\ind\Sigma_1$
   such that $(M,A),(M,B)\models\bd\Sigma^0_1$,
    and $B$ is arithmetically but not $\Delta^0_1$-definable in $(M,A)$.

A coloring $C$ of pairs of numbers $(s, x)$ is {\it stable} if $\text{lim}_s \ C(s, x)$ exists for all $x$. The principle of Stable Ramsey's Theorem for Pairs ($\mathsf{SRT}^2_2$) states that every stable two-coloring of pairs of numbers has an infinite homogeneous set, i.e.~one in which every pair of numbers has the same color.

\begin{question}
If $\M=(M, \X)$ is a $\bd\Sigma^0_2$-model and $A\in\X$, does every $\Delta^0_1(A)$-instance of $\mathsf{SRT}^2_2$ have an arithmetically definable solution relative to $A$ preserving $\bd\Sigma^0_2$?
\end{question}

A $\bd\Sigma^0_2$-model $\M_0$ that  gives a positive answer to the above question is given in Chong, Slaman and Yang \cite{CSY2014}. In that model, the set of standard natural numbers  is a $\Sigma_2$-cut and every instance of $\mathsf{SRT}^2_2$ has a solution that is low.  $\M_0$ has the special feature  that every definable set of natural numbers is coded in the model, a property not shared by every $\bd\Sigma^0_2$-model.

  \begin{question}
If $\M\models\RCA+\bd\Sigma^0_2+\neg\ind\Sigma^0_2(\emptyset')$, is there a set of minimal degree which is low?
\end{question}
The challenges encountered in implementing a priority construction for a set of minimal  $\aleph^L_\omega$- or $\aleph^L_{\omega_1}$-degree are mirrored in  $\M$.  A solution to one situation could shed light on  the other.

\end{document}